\documentstyle[11pt]{article}
 \newlength{\standardunitlength}
 \newtheorem{lemma}{Lemma}
\newtheorem{theorem}{Theorem} 
\newenvironment{proof}{\noindent {\sc Proof:}}{$\Box$ \vspace{2 ex}}

\begin{document}

\begin{center} Semisimple orbits of Lie algebras and card-shuffling
measures on Coxeter groups \end{center}

\begin{center}
By Jason Fulman
\end{center}

\begin{center}
Dartmouth College
\end{center}

\begin{center}
Department of Mathematics
\end{center}

\begin{center}
6188 Bradley Hall
\end{center}

\begin{center}
Hanover, NH 03755, USA
\end{center}

\begin{center}
email:fulman@math.stanford.edu
\end{center}

1991 AMS Primary Subject Classifications: 20G40, 20F55

\newpage
Proposed running head: Semisimple orbits of Lie algebras

\newpage \begin{abstract} Random walk on the chambers of hyperplanes
arrangements is used to define a family of card shuffling measures
$H_{W,x}$ for a finite Coxeter group $W$ and real $x \neq 0$. By
algebraic group theory, there is a map from the semisimple orbits of
the adjoint action of a finite group of Lie type on its Lie algebra to
the conjugacy classes of the Weyl group. Choosing such a semisimple
orbit uniformly at random thereby induces a probability measure on the
conjugacy classes of the Weyl group. For types $A$, $B$, and the
identity conjugacy class of $W$ for all types, it is proved that for
$q$ very good, this measure on conjugacy classes is equal to the
measure arising from $H_{W,q}$.  \end{abstract}


\section{Introduction}

The first part of this paper defines signed measures $H_{W,x}$ for a
finite Coxeter group $W$ and real $x \neq 0$. By a signed measure is
meant an element of the group algebra of $W$ whose coefficients sum to
one. For type $A$ these measures were discovered by Bayer and Diaconis
\cite{BD} in their analysis of riffle shuffling. Their work was
extended to type $B$ in \cite{BB}. (It is amusing to note as in
\cite{BD} that for tarot cards, which often have up/down directions,
type $B$ shuffling is a better model than type $A$ shuffling). For
types $A$ and $B$ these measures also arise in explicit versions of
the Poincar\'e-Birkhoff-Witt theorem \cite{BergeronWolfgang} and in
splittings of Hochschild homology \cite{H}. Section 3.8 of \cite{SSt}
describes the type $A$ measure in the language of Hopf algebras.

The key tool in defining the measures $H_{W,x}$ will be the theory of
random walks on the chambers of hyperplane arrangements, as initiated
in \cite{BHR} and developed in \cite{BrD}. As noted in \cite{F2} (a
follow-up to this paper), the measures defined here generalize to any
real hyperplane arrangement. The point of Section \ref{definition} is
to focus on the case of arrangements coming from finite Coxeter
groups.

Section \ref{semisimple} connects the measures $H_{W,x}$ with the
finite groups of Lie type. Our interest in this direction arose from a
result in \cite{DMP} which connected type $A$ card shuffling with the
factorization of random polynomials. The result in \cite{DMP} was not
expressed in Lie theoretic terms and is not in general related to the
measures $H_{W,x}$. Rather than studying semisimple conjugacy classes,
we study semisimple orbits of the adjoint action of a finite group of
Lie type on its Lie algebra, and with necessary restrictions on the
characteristic. As discussed in \cite{F4}, a construction using the
affine Weyl group is needed to study semisimple conjugacy
classes. This leads to new models of card shuffling.

As mentioned in the abstract, there is a natural map $\Phi$ from the
semisimple orbits of the adjoint action of a finite group of Lie type
on its Lie algebra to the conjugacy classes of the Weyl
group. Choosing such a semisimple orbit uniformly at random gives a
probability measure on the conjugacy classes of the Weyl group. For
$q$ very good, we show that in some cases this measure on conjugacy
classes is equal to the measure arising from $H_{W,q}$. A long term
goal is to refine this map $\Phi$ so that it associates to each
semisimple orbit an {\it element} of $W$. Furthermore choosing an
orbit at random and applying the refined map should give the measures
$H_{W,q}$. In Section \ref{refine} of this paper we indicate how to do
this unnaturally for types $A$ and $B$. A refinement of $\Phi$ which
is both natural and general remains elusive, but could have important
applications in algebraic number theory. The paper \cite{F4} treats
analogous issues for semisimple conjugacy classes.

\section{Definition and Properties of $H_{W,x}$} \label{definition}

To begin we review work of Bidigare, Hanlon, and Rockmore
\cite{BHR}. Let $\cal{A}$ $= \{H_i : i \in I\}$ be a central
hyperplane arrangement (i.e. $\cap_{i \in I} H_i = 0$) for a real
vector space $V$. Let $\gamma$ be a vector in the complement of
$\cal{A}$. Every $H_i$ partitions $V$ into three pieces: $H_i^0=H_i$,
the open half space $H_i^+$ of $V$ containing $\gamma$, and the open
half space $H_i^-$ of $V$ not containing $\gamma$. The faces of
$\cal{A}$ are defined as the non-empty intersections of the form \[
\cap_{i \in I} H_i^{\epsilon_i} \] where $\epsilon_i \in
\{0,-,+\}$. Equivalently, $\cal{A}$ cuts $V$ into regions called
chambers and the faces are the faces of these chambers viewed as
polyhedra.

	A random process (henceforth called the BHR walk) on chambers
is then defined as follows. Assign weights $v(F)$ to the faces of
$\cal{A}$ in such a way that $v_F \geq 0$ for all $F$ and $\sum_F
v(F)=1$. Pick a starting chamber $C_0$. At step $i$, pick a face $F_i$
with chance of face $F$ equal to $v(F)$ and define $C_i$ to be the
chamber adjacent to $F_i$ which is closest to $C_{i-1}$ (separated
from $C_{i-1}$ by the fewest number of hyperplanes.) Such a chamber
always exists.

	To give our definition of $H_{W,x}$, some additional notation
is needed. Let $L$ be the set of intersections of the hyperplanes in
$\cal{A}$, taking $V \in L$. Partially order $L$ by reverse
inclusion. (This lattice is not the same as the face lattice). Recall
that the Moebius function $\mu$ is defined by $\mu(X,X)=1$ and
$\sum_{X \leq Z \leq Y} \mu(Z,Y)=0$ if $X<Y$ and $\mu(X,Y)=0$
otherwise. The characteristic polynomial of $L$ is defined as

	\[ \chi(L,x) = \sum_{X \in L} \mu(V,X) x^{dim(X)}. \] Let $\Pi$ be a base of the
 positive roots of $W$. For $J \subseteq \Pi$, let $Fix(W_J)$ denote
the fixed space of the parabolic subgroup $W_J$ in its action on
$V$. Let $L^{Fix(W_J)}$ denote the restricted poset $\{Y \in
L(\cal{A})$$ |Y \geq Fix(W_J)\}$. Define $Des(w)$ to be the simple
positive roots which $w$ maps to negative roots (also known as the
descent set of $w$) and let $d(w)=|Des(w)|$. Let $N_W(W_K)$ be the
normalizer of $W_K$ in $W$ and let $\lambda(K)$ be the subsets of
$\Pi$ equivalent to $K$ under the action of $W$.

{\bf Definition:} For $W$ a finite Coxeter group and $x \neq 0$,
define $H_{W,x}(w)$ to be \[ \sum_{K \subseteq \Pi-Des(w)}
 \frac{|W_K| \chi(L^{Fix(W_K)},x)}{x^n|N_W(W_K)||\lambda(K)|}.\]

We remark that the paper \cite{BHR} had a hyperplane definition for type $A$
shuffling, but not using group theoretically defined face weights.

To give a feeling for these measures and for later use, we recall
formulas for types $A$ and $B$ (obtained using descent algebras and
also arising from the above definition).

\begin{itemize}

\item (\cite{BD}) \[ H_{S_n,x}(w) = \frac{{x+n-1-d(w) \choose
n}}{x^n}. \] Physically, the inverse of this measure is obtained by
cutting at card $k$ with probability $\frac{{n \choose k}}{2^n}$, then
doing a uniformly chosen random interleaving of the piles. The papers
\cite{DMP} and \cite{F} investigate the cycle structure and inversion
structure of a random permutation chosen from $H_{S_n,x}$.

\item (\cite{BB}) \[ H_{B_n,x}(w) = \frac{(x+2n-1-2d(w))(x+2n-3-2d(w))
\cdots (x+1-2d(w))}{x^n n!}.\] The inverse of this measure also has a
physical description if $x$ is odd, verified for $x=3$ in
\cite{BB}. One cuts multinomially into an odd number of piles, flips
over the even numbered piles, and then does a random
interleaving. This is different from the type $B$ notion in \cite{BD},
which cuts into two piles. However these two types of shuffles can be
placed in a unified setting, using the affine Weyl group \cite{F4}. In
future work we hope to study physical models of the shuffles $H_{W,x}$
for other finite Coxeter groups, viewed as permutation groups.

\end{itemize}

Next we comment on some properties of the measures $H_{W,x}$.

\begin{itemize}

\item (\cite{F2}) For types $A,B,C,H_3$ and rank 2 groups (but not for
all types as is explained in below), the measures $H_{W,x}$ convolve
in the sense that

\[ \left( \sum_{w \in W} H_{W,x}(w)w \right) \left( \sum_{w \in W}
H_{W,y}(w)w \right) = \sum_{w \in W} H_{W,xy}(w) w. \] Thus $n$
$x$-shuffles is the same as an $x^n$ shuffle. Observe also that in the
$x \rightarrow \infty$ limit the measures $H_{W,x}$ become the uniform
distribution. The eigenvalues of an $x$-shuffle viewed as a Markov
chain are $\frac{1}{x^i}$ ,$i=0,\cdots,n-1$ with various
multiplicities.

\item The Coxeter complex of $W$ has as faces the left cosets $wW_K$
and as chambers the elements of $W$. Consider the BHR walks on the
chambers of the Coxeter complex with face wieghts \[v(wW_K) =
\frac{|W_K| \chi(L^{Fix(W_K)},x)} {x^n |N_W(W_K)||\lambda(K)|}.\] When
these weights are non-negative, $H_{W,x}(w)$ can be interpreted as the
probability of moving from the identity chamber to $w$. Equations from
page 282 of \cite{OS} imply that, $v(wW_K)$ can be rewritten as
$(-1)^{n-|K|} \frac{\chi(L^{Fix(W_K)},x)}{x^n
\chi(L^{Fix(W_K)},-1)}$. As observed in \cite{F2}, this leads to a
notion of card shuffling for any real hyperplane arrangement or
oriented matroid. The Coxeter case gives rise to the factorization \[
\chi(L^{Fix(W_K)},x) = \prod_{i=1}^{dim(Fix(W_K))} (x-b_i^K) \] from
\cite{OS} where the $b_i^K$ are integers called coexponents. From the
results and tables in \cite{OS}, all $b_i^K$ are less than or equal to
the maximum exponent of $W$. From the table of bad primes for
crystallographic types on page 28 of \cite{C}, the bad primes are
precisely the primes less than the maximum exponent of $W$ which are
not equal to exponents of $W$. (Equivalently, a prime is good if it
divides no coefficient of any root expressed as a linear combination
of simple roots.) Thus $H_{W,q}(w) \geq 0$ if $W$ is crystallographic
and $q$ is a good prime, because then every face weight is
non-negative. This may be regarded as evidence in favor of the
the statement in Problem 1 in Section \ref{semisimple}.

\item Orlik and Solomon \cite{OS} have calculated and tabulated
$\chi(L^{Fix(W_K)},x)$ for all types. By the previous remark, this
gives a simple and unified method for computing the measure
$H_{W,x}$. Applied to $W$ of type $H_3$, one concludes that

\[ H_{H_3,x}(w) =             \left\{ \begin{array}{ll}
																																				\frac{(x+9)(x+5)(x+1)}{120x^3} &
\mbox{if $d(w)=0$}\\
																																				\frac{(x+5)(x+1)(x-1)}{120x^3} & \mbox{if
$d(w)=1$}\\
										                          \frac{(x+1)(x-1)(x-5)}{120x^3} & \mbox{if
$d(w)=2$}\\
																																				\frac{(x-1)(x-5)(x-9)}{120x^3}		& \mbox{if
$d(w)=3$}.
																																				\end{array}
			\right.			 \]

	This formula, together with the formulas for $H_{W,x}$ for $W$
of types $A,B$ which appeared earlier in this paper, suggest that
$H_{W,x}$ satisfies the following factorization and reciprocity
properties:

\begin{enumerate}

\item $H_{W,x}(w)$ splits into linear factors as a function of $x$.

\item $H_{W,x}(w) = H_{W,-x}(ww_0)$ where $w_0$ is the longest
element of $W$.

\end{enumerate}

	In fact neither of these properties holds. This is evident
from the following formula for $H_{H_4,x}$ which is obtained by using
tables of Orlik and Solomon as just described.

\[ H_{H_4,x}(w) =             \left\{ \begin{array}{ll}
																																				\frac{(x+29)(x+19)(x+11)(x+1)}{14400x^4} &
\mbox{if $d(w)=0$}\\
																																				\frac{(x+1)(x-1)(x^2+30x+149)}{14400x^4} &
\mbox{if $Des(w)=\{\alpha_1\}$ or $Des(w)=\{\alpha_2\}$}\\
										                          \frac{(x+1)(x-1)(x^2+30x+269)}{14400x^4} &
\mbox{if
$Des(w)=\{\alpha_3\}$ or $Des(w)=\{\alpha_4\}$}\\
																																				\frac{(x+11)(x+1)(x-1)(x-11)}{14400x^4} &
\mbox{if $d(w)=2$ and $Des(w) \neq \{\alpha_3,\alpha_4\}$}\\
																																				\frac{(x+1)^2(x-1)^2}{14400x^4} &
\mbox{if $Des(w)=\{\alpha_3,\alpha_4\}$}\\
							                             \frac{(x+1)(x-1)(x-11)(x-19)}{14400x^4} &
\mbox{if $d(w)=3$}\\

                                    \frac{(x-1)(x-11)(x-19)(x-29)}{14400x^4}		&
\mbox{if
$d(w)=4$}.
																																				\end{array}
			\right.  \] Incidentally this remark shows
that $H_{H_4,x}$ does not convolve. For $H_{W,-1}$ places all mass on
the longest element $w_0$, so the convolution property would imply
that $H_{H_4,-x}(w)=H_{H_4,-x}(ww_0)$. Since $w$ and $ww_0$ have
complementary descent sets, this equality does not hold for $w$ with
$Des(w)=\{\alpha_3,\alpha_4\}$. The same argument disproves the
convolution property in many cases.

\end{itemize}

Let $id$ be the identity element of $W$ and $w_0$ the longest element of
$W$. Theorem \ref{longshort} calculates the values of the measure $H_{W,x}$ on
these elements.

\begin{theorem} \label{longshort} Let $m_1,\cdots,m_r$ be the exponents of $W$.
Then
\begin{eqnarray*}
H_{W,x}(w_0) & = & \frac{\prod_{i=1}^r (x-m_i)}{x^r |W|}\\
H_{W,x}(id) & = & \frac{\prod_{i=1}^r (x+m_i)}{x^r |W|}.
\end{eqnarray*}
\end{theorem}

\begin{proof}
The first assertion is easier. In fact,

\[H_{W,x}(w_0)  =  \frac{\chi(L,x)}{x^r |W|} = \frac{\prod_{i=1}^r (x-m_i)}{x^r
|W|}.\]The first equality is from the definition of $H_{W,x}$ and the second
equality is a well known factorization of the characteristic polynomial of
$L$ (e.g. \cite{OS}).

For the second assertion, additional concepts are needed. Let $L$ be
the lattice in $V$ generated by $\check{\Phi}$ and let \[ \hat{L} =
\{v \in V | <v,\alpha> \in Z \ for \ all \ \alpha \in \Phi \}. \] Let
$f = [\hat{L}:L]$ be the index $L$ in $\hat{L}$. Let $\Pi =
\{\alpha_i\} \subset \Phi^+$ be a set of simple roots contained in a
set of positive roots and let $\theta$ be the highest root in
$\Phi^+$. For convenience set $\alpha_0 = - \theta$. Let $\tilde{\Pi}=
\Pi \cup \{ \alpha_0 \}$. Define coefficients $c_{\alpha}$ of $\theta$
with respect to $\tilde{\Pi}$ by the equations $\sum_{\alpha \in
\tilde{\Pi}} c_{\alpha} \alpha = 0$ and $c_{\alpha_0}=1$. For $S \neq
\tilde{\Pi}$ a proper subset of $\tilde{\Pi}$, define as in \cite{So}
$p(S,x)$ to be the number of solutions ${\bf y}$ in strictly positive
integers to the equation

\[ \sum_{\alpha \in \tilde{\Pi}-S} c_{\alpha} y_{\alpha} = x. \]

In the equations which follow $W_{K_1},\cdots,W_{K_l}$ with
$K_1,\cdots,K_l \subseteq \Pi$ are representatives for the parabolic
subgroups of $W$ under conjugation. In \cite{OS} it is proved that
$|\lambda(K)|$ is the number of $J \subseteq \Pi$ such that $W_J$ is
conjugate to $W_K$. We also make use of the fact \cite{So2} that if
$x$ is relatively prime to all $c_{\alpha}$, then for any $S \subset
\tilde{\Pi}, S \neq \tilde{\Pi}$, if $p(S,x)$ is non-zero then $W_S$
is conjugate to one of $W_{K_1},\cdots,W_{K_l}.$ We denote conjugacy
of parabolic subgroups by the symbol $\sim$. One concludes that for
infinitely many (and hence all) non-zero $x$,

\begin{eqnarray*}
H_{W,x}(id) & = & \sum_{K \subseteq \Pi} \frac{|W_K|
\chi(L^{Fix(W_K)},x)}{x^r|N_W(W_K)||\lambda(K)|}\\
& = & \sum_{i=1}^l \frac{|W_{K_i}| \chi(L^{Fix(W_{K_i})},x)}{x^r|N_W(W_{K_i})|}\\
&=& \frac{1}{x^rf} \sum_{i=1}^l \frac{f |W_{K_i}|
\chi(L^{Fix(W_{K_i})},x)}{|N_W(W_{K_i})|}\\
& = & \frac{1}{x^rf} \sum_{i=1}^l \sum_{S \subseteq \tilde{\Pi},S \neq \tilde{\Pi}
\atop W_S \sim W_{K_i}} p(S,x)\\
& = & \frac{1}{x^rf} \sum_{S \subset \tilde{\Pi} \atop S \neq \tilde{\Pi}}
p(S,x)\\
& = & \frac{1}{x^r|W|} \prod_{i=1}^r (x+m_i).
\end{eqnarray*} The fourth and sixth equalities are results of \cite{So}.
\end{proof}

\section{Semisimple Orbits of Lie Algebras} \label{semisimple}

	This section connects the signed measures $H_{W,x}$ with
semisimple orbits of the adjoint action of finite groups of Lie type
on their Lie algebras.

	 Let $G$ be a connected semisimple group defined over a finite
field of $q$ elements. Suppose also that $G$ is simply connected. Let
$\cal G$ be the Lie algebra of $G$. Let $F$ denote both a Frobenius
automorphism of $G$ and the corresponding Frobenius automorphism of
$\cal G$. Suppose that $G$ is $F$-split.  Since the derived group of
$G$ is simply connected (the derived group of a simply connected group
is itself), a theorem of Springer and Steinberg \cite{SS} implies that
the centralizers of semisimple elements of $\cal G$ are connected.
Let $r$ be the rank of $G$.

	Now we define a map $\Phi$ (studied in \cite{L} in somewhat
greater generality) from the $F$-rational semisimple orbits $c$ of
$\cal G$ to $W$, the Weyl group of $G$.  Pick $x \in {\cal G}^F \cap
c$. Since the centralizers of semisimple elements of $\cal G$ are
connected, $x$ is determined up to conjugacy in $G^F$ and $C_G(x)$,
the centralizer in $G$ of $x$, is determined up to $G^F$
conjugacy. Let $T$ be a maximally split maximal torus in
$C_G(x)$. Then $T$ is an $F$-stable maximal torus of $G$, determined
up to $G^F$ conjugacy. By Proposition 3.3.3 of \cite{C}, the $G^F$
conjugacy classes of $F$-stable maximal tori of $G$ are in bijection
with conjugacy classes of $W$. Define $\Phi(c)$ to be the
corresponding conjugacy class of $W$.

	For example, in type $A_{n-1}$ the semisimple orbits $c$ of
$sl(n,q)$ correspond to monic degree $n$ polynomials $f(c)$ whose
coefficient of $x^{n-1}$ vanishes.  Such a polynomial factors as
$\prod_i f_i^{a_i}$ where the $f_i$ are irreducible over
$F_q$. Letting $d_i$ be the degree of $f_i$, $\Phi(c)$ is the
conjugacy class of $S_n$ corresponding to the partition $(d_i^{a_i})$.

	As is standard in Lie theory (e.g. \cite{Fl}), call a prime
$p$ very good if it divides no coefficient of any root expressed as a
linear combination of simple roots and is relatively prime to the
index of connection (the index of the coroot lattice in the weight
lattice). For example in type $A$ the very good primes are those not
dividing $n$.

\vspace{.5mm}

{\bf Problem 1:} When is the following statement true? "Let $G$ be as
above, and suppose that the characteristic is a prime which is very
good for $G$. Choose $c$ among the $q^r$ $F$-rational semisimple
orbits of $\cal G$ uniformly at random. Then for all conjugacy classes
$C$ of $W$,

\[ Prob(\Phi(c) = C) = \sum_{w \in C} H_{W,q}(w)." \]

	Recall from the end of Section \ref{definition} that under the
conditions of Problem 1, $H_{W,q}(w) \geq 0$ for all $w \in W$. This
may be taken as evidence that the statement in Problem 1 is
correct. Theorems \ref{conj2id}, \ref{conj2sym}, and \ref{conj2hyp}
provide further evidence. In cases where the convolution property of
$W$ does not hold, we have doubts as to whether the statement in
Problem 1 is always true. Nevertheless, at present we have no examples
to the contrary (though type $D_4$ would be a natural first place to
look).

\begin{theorem} \label{conj2id} The statement in Problem 1 holds for
$G$ of all types (i.e. $A-D,E_{6-8},F_4,G_2$) when $C$ is the identity
conjugacy class of $W$.  \end{theorem}

\begin{proof}
	Corollary 3.4 of \cite{Fl} (see also Proposition 5.9 of
\cite{L}) states that for $q$ very good, the number of $F$-rational
semisimple orbits $c$ of $\cal G$ which satisfy $\Phi(c)=id$ is equal
to \[ \prod_{i=1}^r \frac{q+m_i}{1+m_i} \] where $r$ is the rank of
$G$ and $m_i$ are the exponents of $W$. Since there are a total of
$q^r$ $F$-rational semisimple orbits of $\cal G$, and because
$|W|=\prod_{i=1}^r (1+m_i)$, \[ Prob( \Phi(c) = id) =
\frac{\prod_{i=1}^r (q+m_i)}{q^r|W|}. \] The proposition now follows
from Theorem \ref{longshort}.  \end{proof}

\begin{theorem} \label{conj2sym} The statement of Problem 1 holds for
$G$ of type $A$, for all conjugacy classes $C$ of the symmetric group
$S_n$.  \end{theorem}

\begin{proof}
	Note that a monic, degree $n$ polynomial $f$ with coefficients
in $F_q$ defines a partition of $n$, and hence a conjugacy class of
$S_n$, by its factorization into irreducibles. To be precise, if $f$
factors as $\prod_i f_i^{a_i}$ where the $f_i$ are irreducible of
degree $d_i$, then $(d_i^{a_i})$ is a partition of $n$.  If the
coefficient of $x^{n-1}$ in $f$ vanishes, then $f$ represents an
$F$-rational semisimple orbit $c$ of $sl(n,q)$, and the conjugacy
class of $S_n$ corresponding to the partition $(d_i^{a_i})$ is equal
to $\Phi(c)$.

	In \cite{DMP} it is shown that if $f$ is uniformly chosen
among all monic, degree $n$ polynomials with coefficients in $F_q$,
then the measure on the conjugacy classes of $S_n$ induced by the
factorization of $f$ is equal to the measure induced by
$H_{S_n,q}$. Thus, to prove the theorem, it suffices to show that the
random partition associated to a uniformly chosen monic, degree $n$
polynomial over $F_q$ has the same distribution as the random
partition associated to a uniformly chosen monic, degree $n$
polynomial over $F_q$ with vanishing coefficient of $x^{n-1}$. Since
the characteristic $p$ is assumed to be very good, $p$ does not divide
$n$. Thus for a suitable choice of $k$, the change of variables $x
\rightarrow x+k$ gives rise to a bijection between monic, degree $n$
polynomials with coefficient of $x^{n-1}$ equal to $b_1$ and monic,
degree $n$ polynomials with coefficient of $x^{n-1}$ equal to $b_2$,
for any $b_1$ and $b_2$. Since this bijection preserves the partition
associated to a polynomial, the theorem is proved.  \end{proof}

	Theorem \ref{conj2hyp} will confirm the statement of Problem 1
for all $G$ of type $B$. The proof will use the following
combinatorial objects introduced in \cite{R}. As Lemma \ref{equiv}
will show, these objects have interpretations in terms of
polynomials. Let a ${\bf Z}$-word of length $m$ be a vector
$(a_1,\cdots,a_m) \in {\bf Z}^m$. For such a word define $max(a)=
max(|a_i|)_{i=1}^m$. The cyclic group $C_{2m}$ acts on ${\bf Z}$-words
of length $m$ by having a generator $g$ act as $g(a_1,\cdots,a_m) =
(a_2,\cdots,a_m,-a_1)$. Call a fixed-point free orbit $P$ of this
action a primitive twisted necklace of size $m$. The group $Z_2 \times
C_m$ acts on ${\bf Z}$-words of length $m$ by having the generator $r$
of $C_m$ act as a cyclic shift
$r(a_1,\cdots,a_m)=(a_2,\cdots,a_m,a_1)$ and having the generator $v$
of $Z_2$ act by $v(a_1,\cdots,a_m)=(-a_1,\cdots,-a_m)$. Call an orbit
$D$ of this action a primitive blinking necklace of size $m$ if its
$C_m$ action is free (though its $Z_2 \times C_m$ action need not be).
Let a signed ornament $o$ be a set of primitive twisted necklaces and
a multiset of primitive blinking necklaces. Say that $o$ has type
$(\vec{\lambda},\vec{\mu}) =
((\lambda_1,\lambda_2,\cdots),(\mu_1,\mu_2,\cdots))$ if it consists of
$\lambda_m$ primitive blinking neclaces of size $m$ and $\mu_m$
primitive twisted necklaces of size $m$. Also define the size of $o$
to be the sum of the sizes of the primitive twisted and blinking
necklaces which make up $o$, and define $max(o)$ to be the maximum of
$max(D)$ and $max(P)$ for the primitive twisted and blinking necklaces
which make up $o$.

\begin{lemma} \label{equiv} Primitive twisted necklaces $P$ of size
$m$ and with $max(P) \leq \frac{q-1}{2}$ correspond to irreducible
polynomials $f(z)$ over $F_q$ of degree $2m$ satisfying
$f(z)=f(-z)$. Primitive blinking necklaces $D$ of size $m$ and with
$max(D) \leq \frac{q-1}{2}$ correspond to products $f(z)f(-z)$ with
$f(z),f(-z)$ a pair of irreducible polynomial of degree $m$ over
$F_q$. Signed ornaments given as sets of such $P$'s and multisets of
such $D$'s correspond to polynomials of degree $2m$ over $F_q$
satisfying $f(z)=f(-z)$.  \end{lemma}

\begin{proof} For the first assertion, let $F_{q^{2m}}$ be the degree $2m$
extension of $F_q$. Choose $\alpha$ in $F_{q^{2m}}$ such that
$\{\alpha^{q^i}:1 \leq i \leq 2m\}$ is a basis over $F_q$. (Such a basis is
called a normal basis and is known to exist). Let $f(z)$ be an irreducible
polynomial of degree $2m$ satisfying $f(z)=f(-z)$. Let $\beta$ be one of its
roots in $F_{q^{2m}}$. Writing $\beta = \sum_{i=1}^{2m} c_i \alpha^{q^i}$,
define a vector $(c_1,\cdots,c_m)$ associated to $\beta$. Since the automorphism
of $F_{q^{2m}}$ defined by $\alpha \rightarrow \alpha^{q^m}$ is its unique
automorphism of order two, it follows that $\beta^q$ is assigned the vector
$(c_2,\cdots,c_m,-c_1)$. Thus the action of the Frobenius map
$x \rightarrow x^q$ corresponds to the action of $Z_2 \times C_m$ on the vector
$(c_1,\cdots,c_m)$, and irreducible polynomials correspond to primitive orbits.

For the second assertion, choose $\alpha$ in $F_{q^{m}}$ such
that $\{\alpha^{q^i}:1 \leq i \leq m\}$ is a basis over $F_q$. Let $f(z)$ be an
irreducible polynomial of degree $m$. Let $\beta$ be one
of its roots in $F_{q^{m}}$. Writing $\beta = \sum_{i=1}^{m} c_i
\alpha^{q^i}$, define a vector $(c_1,\cdots,c_m)$ associated to $\beta$. The
$C_m$ action on this vector is free because $f(z)$ is irreducible. The $Z_2$
action sends $f(z)$ to $f(-z)$.

For the final assertion, note that a polynomial $f(z)$ satisfying
$f(z)=f(-z)$ can be factored uniquely as a product \[
\prod_{\{\phi_j(z), \phi_j(-z)\}} [\phi_j(z) \phi_j(-z)]^{r_{\phi_j}}
\prod_{\phi_j : \phi_j(z) = \phi_j(-z)} \phi_j(z)^{s_{\phi_j}} \]
where the $\phi_j$ are monic irreducible polynomials and $s_{\phi_j}
\in \{0,1\}$.
\end{proof}

Theorem \ref{conj2hyp} proves the statement of Problem 1 for type
$B$.

\begin{theorem} \label{conj2hyp} The statement of Problem 1 holds for
$G$ of type $B$, for all conjugacy classes $C$ of the hyperoctahedral
group $B_n$.  \end{theorem}

\begin{proof}
	Note that because $2$ is a bad prime for type $B$, it can be
assumed that the characteristic is odd. Recall that the type of a
signed ornament is parameterized by pairs of vectors
$(\vec{\lambda},\vec{\mu})$, where $\lambda_i$ is the number of
primitive blinking necklaces of size $i$ and $\mu_i$ is the number of
primitive twisted necklaces of size $i$. From the theory of wreath
products the conjugacy classes of the hyperoctahedral group $B_n$ are
also parameterized by pairs of vectors $(\vec{\lambda},\vec{\mu})$,
where $\lambda_i(w)$ and $\mu_i(w)$ are the number of positive and
negative cycles of $w \in B_n$ respectively.

	The first step of the proof will be to show that the measure
induced on pairs $(\vec{\lambda},\vec{\mu})$ by choosing a random
signed ornament $o$ of size $n$ satisfying $max(o) \leq \frac{q-1}{2}$
is equal to the measure induced on pairs $(\vec{\lambda},\vec{\mu})$
by choosing $w \in B_n$ according to the measure $H_{B_n,q}$ and then
looking at its conjugacy class. From the definition of descents given
in Section \ref{definition}, it is easy to see that if one introduces
the following linear order $\Lambda$ on the set of non-zero integers:
\[ +1 <_{\Lambda} +2 <_{\Lambda} \cdots +n<_{\Lambda} \cdots
<_{\Lambda} -n <_{\Lambda} \cdots <_{\Lambda} -2 <_{\Lambda} -1 \]
then $d(w)$, the number of descents of $w \in B_n$, can be defined as
$|\{i: 1 \leq i \leq n: w(i) <_{\Lambda} w(i+1)\}|$. Here $w(n+1)=n+1$
by convention.

	It is proved in \cite{R} that there is a bijection between signed
ornaments
$o$ of size $n$ satisfying $max(o) \leq \frac{q-1}{2}$ and pairs $(w,\vec{s})$
where
$w
\in B_n$ and $\vec{s}=(s_1,\cdots,s_n) \in {\bf N}^n$ satisifies $\frac{q-1}{2}
\geq s_1 \geq \cdots \geq s_n \geq 0$ and $s_i > s_{i+1}$ when $w(i) <_{\Lambda}
w(i+1)$ (i.e. when $w$ has a descent at position $i$). Further, he shows that
the type of $o$ is equal to the conjugacy class vector of $w$. It is easy to see
that if $w$ has $d(w)$ descents, then the number of
$\vec{s}$ such that $\frac{q-1}{2}
\geq s_1 \geq \cdots \geq s_n \geq 0$ and $s_i > s_{i+1}$ when $w(i) <_{\Lambda}
w(i+1)$ is equal to

\[ {\frac{q-1}{2}+n-d(w) \choose n} = \frac{(q+1-2d(\pi)) \cdots
(q+2n-1-2d(\pi))}{2^n n!}. \]

	Lemma \ref{equiv} implies that there are $q^n$ signed ornaments $f$ of size
$n$ satisfying $max(f) \leq \frac{q-1}{2}$. Thus choosing a random signed
ornament induces a measure on $w \in B_n$ with mass on $w$ equal to

\[ \frac{(q+1-2d(\pi)) \cdots (q+2n-1-2d(\pi))}{q^n |B_n|}. \] By the
remarks in Section \ref{definition}, this is exactly the mass on $w$
under the measure $H_{B_n,q}$. Since in Reiner's bijection the type of
$o$ is equal to the conjugacy class vector of $w$, we have proved that
the measure on conjugacy classes $(\vec{\lambda},\vec{\mu})$ of $B_n$
induced by choosing $w$ according to $H_{B_n,q}$ is equal to the
measure on conjugacy classes $(\vec{\lambda},\vec{\mu})$ of $B_n$
induced by choosing a signed ornament uniformly at random and taking
its type.

	The second step in the proof is to show that if $f$ is chosen
uniformly among the $q^n$ semisimple orbits of $Spin(2n+1,q)$ on its
Lie algebra, then the chance that $\Phi(f)$ is the conjugacy class
$(\vec{\lambda},\vec{\mu})$ of $B_n$ is equal to the chance that a
signed ornament chosen randomly among the $q^n$ signed ornaments $o$
of size $n$ satisfying $max(o) \leq \frac{q-1}{2}$ has type
$(\vec{\lambda},\vec{\mu})$. It is well known that the semisimple
orbits of $Spin(2n+1,q)$ on its Lie algebra correspond to monic,
degree $2n$ polynomials $f$ satisfying $f(z)=f(-z)$. From Section 2 of
\cite{SS} and Section 3 of \cite{C2}, one sees that $\Phi(f)$ can be
described as follows. Factor $f$ uniquely into irreducibles as

\[ \prod_{\{\phi_j(z), \phi_j(-z)\}} [\phi_j(z)
\phi_j(-z)]^{r_{\phi_j}} \prod_{\phi_j : \phi_j(z) = \phi_j(-z)}
\phi_j(z)^{s_{\phi_j}} \] where the $\phi_j$ are monic irreducible
polynomials and $s_{\phi_j} \in \{0,1\}$. Then let
$\lambda_i(f)=\sum_{\phi: deg(\phi)=i} r_{\phi}$ and
$\mu_i(f)=\sum_{\phi: deg(\phi)=2i} s_{\phi}$. The result now follows
from Lemma \ref{equiv}.
\end{proof}

	We remark that the statement of Problem 1 would be false if
instead of choosing $c$ uniformly among the $q^r$ $F$-rational
semisimple orbits of ${\cal G}$, $c$ were chosen uniformly among the
$q^r$ semisimple conjugacy classes of $G^F$. For a simple
counterexample, take $G=SL(3,5)$ and $C$ the identity conjugacy class
of $S_3$. There are only five monic polynomials $f$ with coefficients
in $F_5$ which factor into linear terms and satisfy $f(0)=1$. The
analog of the statement of Problem 1 would predict that there are
seven. For analogous, yet combinatorially more intricate developments
for semisimple conjugacy classes, see \cite{F4}.

\section{Refining the Map $\Phi$ to the Weyl group} \label{refine}

	As noted in the introduction, one long-term goal is to find a
canonical way to associate to an $F$-rational semisimple orbit $c$ of
${\cal G}$ an {\it element} $w$ of $W$. The conjugacy class of $w$
should equal $\Phi(c)$ and choosing $c$ uniformly at random should
induce the measure $H_{W,q}$ on $W$.

	To see why such a result may be interesting, at least in type
$A$, consider a simple algebraic extension of $Q$ with minimal
polynomial $f(x)$. At unramified primes the Frobenius automorphism is
defined up to conjugacy in the Galois group. Viewed as a permutation
of the roots of $f(x)$, the cycle structure of the Frobenius
automorphism is given by the degrees of the irreducible factors of the
modulo $p$ reduction of $f(x)$. This is simply the map $\Phi$ in type
$A$. Some important constructions in algebraic number theory (see
\cite{Gel} for a survey) create generating functions combining this
data over all primes. It is not impossible that a natural refinement
of the Frobenius data will yield new number theoretic constructions.

	Next we indicate a somewhat unnatural way to refine the map
$\Phi$ in types $A$ and $B$. For type $A$, the refinement proceeds in
two steps. Define a necklace on an alphabet to be a sequence of
cyclically arranged letters of the alphabet. A necklace is said to be
primitive if it is not equal to any of its non-trivial cyclic
shifts. For example, the necklace $(a\ a\ b\ b)$ is primitive, but the
necklace $(a\ b\ a\ b)$ is not.

	The first step is to associate to a monic degree $n$
polynomial over $F_p$ a multiset of primitive necklaces on the
alphabet $\{0,1,\cdots,p-1\}$. One way to do this is using the concept
of a normal basis, that is to choose for each $n$ an element
$\alpha_n$ such that its conjugates $\alpha_n^{p^j}$ for
$j=0,\cdots,n-1$ are a basis of $F_{p^n}$ over $F_p$. Then a monic
irreducible degree $i$ polynomial gives a primitive necklace of size
$i$ formed by the coefficients $c_j$ of any one of its roots written
as $\sum c_j \alpha_i^{p^j}$. (It is natural to require that for
$i|n$, the norm of $\alpha_n$ is $\alpha_i$.) This is the preferred
method in the case of semisimple adjoint orbits, because the
involution sending $f(x)$ to $f(-x)$ takes negatives of the necklace
entries.

	A second way to carry out this first step was noticed by
Golomb \cite{Go}. For each $n$, pick an element $\beta_n$ generating
the multiplicative group of the field extension $F_{p^n}$ of $F_p$. A
root of an irreducible polynomial $\phi$ of degree $i$ can be written
$\beta_i^x$. Considering the mod $p$ expansion of $x$ gives a
primitive necklace of size $i$. This is the preferred construction in
the case of semisimple conjugacy classes, because the involution $f(x)
\mapsto \frac{t^{deg(f)} f (\frac{1}{t})} {f(0)}$ on polynomials with
non-zero constant term takes negatives of the necklace entries.

	The next step in the construction is to associate to a
multiset of primitive necklaces on $\{0,\cdots,p-1\}$ a permutation
with cycle structure equal to that of the necklace. A way to do this
was found by Gessel and Reutenauer \cite{G}, and its importance for
card shuffling was recognized in \cite{DMP}. To each entry of a
necklace, first associate the infinite word obtained by reading the
necklace in the clockwise direction. Using the example from \cite{G},
consider the multiset of necklaces \[ (1\ 2)(1\ 2)(2)(2\ 3)(2\ 3\ 2\
3\ 3). \] Then the entry $2$ on the necklace $(2\ 3)$ would give the
word $23232323 \cdots$. One then orders lexicographically the words
obtained (after imposing an arbitrary order on equal necklaces), and
replaces each necklace entry by the lexicographic order of its
associated word. The example would thus yield the permutation \[ (1\
3)(2\ 4)(5)(6\ 9)(7\ 11\ 8\ 12\ 10). \]

	For a $B_n$ analog, the bijection of Gessel should be replaced
by the bijection of Reiner \cite{R} used in the proof of Theorem
\ref{conj2hyp}.

\section{Acknowledgements} The author thanks Persi Diaconis, Dick
Gross, Vic Reiner, and Eric Sommers for helpful discussions.

\end{document}